# PRIME RADICAL AND PRIMARY DECOMPOSITION
# OF IDEALS IN  AN L-SUBRING

**NASEEM AJMAL**

*Department of Mathematics, Zakir Hussain College,*
*(University of Delhi), J.L. Nehru Marg, New Delhi, India.*

**ANAND SWAROOP PRAJAPATI***

*Department of Mathematics, Atma Ram Sanatan Dharma College,*

*(University of Delhi), Dhaula Kuan, New Delhi-110021, (India).*

## ABSTRACT

In this paper we introduce the concept of a prime radical of an ideal of an L-ring $L(\mu, R)$. Among various results pertaining to this concept, we prove here that prime radicals of an ideal $\eta$, its radical $\sqrt{\eta}$, its semiprime radical $S(\eta)$ and its prime radical $P(\eta)$, all coincide. Also we prove that for a primary ideal, its prime radical coincide with its radical. Moreover, we introduce the concept of primary decomposition and reduced primary decomposition of an ideal in an L-ring. We obtain a necessary and sufficient conditions for an ideal of an L-ring to have a primary decomposition. Some more results pertaining to the decomposition of an ideal are established.

* **Correspondence Address :** Anand Swaroop Prajapati, 3764, Motia Bagh, Sarai Phoos, Delhi-110054, India.

## INTRODUCTION

In papers [13,14], we introduced the concept of a maximal ideal of an $L$-ring $L(\mu, R)$. Then in paper [15,16,17], we developed a systematic theory for ideals in an L-ring. We introduced the concepts of a prime ideal, semiprime ideal, radical of an ideal and primary ideal of an $L$-ring. Moreover, several related concepts are formulated such as minimal prime ideal, associated prime ideal and irreducibility of an ideal in an $L$-ring. In paper [18], the notion of semiprime radical of an ideal of an $L$-ring is introduced. All the important ring theoretic analogs of these concepts are established. In another paper [12], the concept of right (left) quotient (or residual) of an ideal $\eta$ by an ideal $\nu$ of an L-ring $\mu$ is introduced and discussed.

With this machinery at our disposal, in this paper, we have further introduced the concept of a prime radical of an ideal of an $L$-ring $L(\mu, R)$. It is proved that the prime radicals of an ideal $\eta$, its radical $\sqrt{\eta}$, its semiprime radical $S(\eta)$ and its prime radical $P(\eta)$, are identical. It is also proved that the prime radical of an ideal of an $L$-ring is always a semiprime ideal. We have also proved that for a primary ideal of an $L$-ring, its radical, semiprime radical and prime radical coincide. We have established that semiprime radical of the prime radical of an ideal of $L$-ring is the prime radical of the ideal.

In Section 3, we have introduced the concepts of primary decomposition (or primary representation) and reduced primary decomposition of an ideal in an $L$-ring. A necessary and sufficient condition for an ideal of an $L$-ring to have a primary decomposition is established. It is also proved that if for each $L$-ring $L(\mu, R)$, every ideal of $\mu$ has a primary decomposition, then each ideal of any subring of $R$ has a primary decomposition in the subring. We have also proved that if an ideal $\eta$ of an $L$-ring $L(\mu, R)$ has a primary decomposition, then every non-empty level subset $\eta_t$ has a primary decomposition in the subring $\mu_t$. Further if, $\eta$ has sup-property and the primary decomposition of some level subset $\eta_t$ is reduced then the primary decomposition of $\eta$ is reduced.

## 1. PRELIMINARIES

In this section, we recall some of the basic definitions and concepts which are used in the sequel. For details we refer to [10,11,13].

Let X be a non-empty set and L be a lattice. By an L-subset of X, we mean a function from X to L. The set of all L-subsets of X is called the L-power set of X and is denoted by $L^X$. For $\mu \in L^X$, the set $\{\mu(x) | x \in X\}$ is called the image of $\mu$, and is denoted by Im$\mu$. For $\mu, \nu \in L^X$, if $\nu(x) \leq \mu(x)$, $\forall \; x \in X$, then we say that $\nu$ is contained in $\mu$ and we write $\nu \subseteq \mu$. If $\nu \subseteq \mu$ and $\nu \neq \mu$, then $\nu$ is said to be properly contained in $\mu$ and we write $\nu \underset{\neq}{\subseteq} \mu$. Throughout the paper, R will denote an ordinary ring and L will denote a lattice unless otherwise specifically mentioned. Also $Z^+$ will denote the set of positive integers and $\phi$ will denote an empty set.

**Definition 1.1.** Let L be a lattice and R be a ring. Let $\mu \in L^R$. Then $\mu$ is called an L-subring of R if

(1)     $\mu(x - y) \geq \mu(x) \wedge \mu(y)$,     $\forall \; x, y \in R$, and

(2)     $\mu(xy) \geq \mu(x) \wedge \mu(y)$,     $\forall \; x, y \in R$.

The set of all L-subrings of R is denoted by L(R). It is obvious that if $\mu$ is an L-subring of R, then $\mu(x) \leq \mu(0)$, $\forall$ $x \in R$. For convenience, we use the notation $L(\mu, R)$ for the L-subring $\mu$ of R and we shall refer to it here as an L-ring $L(\mu, R)$.

**Definition 1.2** Let L be a lattice and $\mu \in L^R$. Then $\mu$ is called L-ideal of R if

(1)   $\mu(x - y) \geq \mu(x) \wedge \mu(y)$,    $\forall$ $x, y \in R$,  and

(2)   $\mu(xy) \geq \mu(x) \vee \mu(y)$ ,    $\forall$ $x, y \in R$.

We denote the set of all L-ideals of R by LI(R). It is obvious that if R has identity 1 and $\mu \in LI(R)$, then $\mu(x) \geq \mu(I)$.

**Definition 1.3.** For $\mu \in L^X$ and $\alpha \in L$, we define level subset $\mu_\alpha$ and strong level subset $\mu_\alpha^>$ of $\mu$ in X, as follows:

$$\mu_\alpha = \{x \in X \mid \mu(x) \geq \alpha\}$$

$$\mu_\alpha^> = \{x \in X \mid \mu(x) > \alpha\}.$$

Obviously $\mu_\alpha^> \subseteq \mu_\alpha$ and for $\alpha \leq \beta$, $\mu_\beta \subseteq \mu_\alpha$ and $\mu_\beta^> \subseteq \mu_\alpha^>$. Also, if $\eta, \theta \in L^X$ with $\eta \subseteq \theta$, then $\eta_\alpha \subseteq \theta_\alpha$, $\forall \alpha \in L$.

**Note :** If $L(\mu, R)$ is an L-ring then each non-empty level subset $\mu_\alpha$ is an ordinary subring of R, called level subring of $\mu$.

**Lemma 1.4.** Let L be a chain and $L(\mu, R)$ be an L-ring. Then each non-empty strong level subset $\mu_r^>$ is a subring R. Also each non-empty level subset $\mu_t$ with $r < t$, is a subring of the subring $\mu_r^>$.

**Definition 1.5 (Definition 3.211 [11]).** Let $\nu \in L^R$ and $\mu \in L(R)$ with $\nu \subseteq \mu$. Then $\nu$ is called an L-ideal of $\mu$ (or in $\mu$) if

(1)     $\nu(x - y) \geq \nu(x) \wedge \nu(y)$,                    $\forall \; x, y \in R$, and

(2)     $\nu(xy) \geq (\mu(x) \wedge \nu(y)) \vee (\nu(x) \wedge \mu(y))$,     $\forall \; x, y \in R$

For convenience, $\nu$ is called an ideal of $\mu$ (or L-ring $L(\mu, R)$).

**Definition 1.6.** Let L be a lattice and $L(\mu, R)$ be an L-ring. If $\nu$ is an L-subring of R with $\nu \subseteq \mu$, then $\nu$ is called a subring of $\mu$ (or L-ring $L(\mu, R)$).

Clearly if $\nu$ is a subring of $\mu$, then $\nu(x^n) \geq \nu(x), \forall n \in Z^+$.

**Theorem 1.7([13]).** Let L be a lattice and R be a ring. Let $L(\mu, R)$ be an L-ring and $\eta \in L^R$ with $\eta \subseteq \mu$. Then $\eta$ is an ideal of $\mu$ if and only if each non-empty level subset $\eta_a$ is an ideal of level subring $\mu_a$.

**Lemma 1.8.** Let L be a lattice and X be a non-empty set. Suppose $\eta, \nu \in L^X$ such that $\eta \subseteq \nu$ with $\eta_t = \nu_t$, $\forall \in \text{Im} \nu$. Then $\eta = \nu$.

**Definition 1.9 [11].** Let L be a complete lattice and $\eta, \nu \in L^R$. Then we define $\eta + \nu$ by

$$(\eta + \nu)(x) = \vee \{\eta(y) \wedge \nu(z) \mid y, z \in R; x = y + z\}$$

Clearly if $\eta$ and $\nu$ are subrings of an L-ring $L(\mu, R)$ with $\eta(0) = \nu(0)$ then $\eta \subseteq \eta + \nu$ and $\nu \subseteq \eta + \nu$.

**Lemma 1.10 ([12]).** Let L be a complete lattice and $L(\mu, R)$ be an L-ring. Let $\eta$ be a subring of $\mu$ then $\eta + \eta = \eta$. Also $\mu + \mu = \mu$.

**Lemma 1.11 ([12]).** Let L be a complete lattice and $L(\mu, R)$ be an L-ring. If $\eta$ and $\nu$ are ideals of $\mu$ with $\eta(0) = \nu(0)$, then $\eta + \nu$ is an ideal of $\mu$ and $\eta \subseteq \eta + \nu$, $\nu \subseteq \eta + \nu$.

**Lemma 1.12 (Theorem 3.2.15 [11]).** Let $L(\mu, R)$ be an L-ring. Then intersection of two ideals of $\mu$ is an ideal of $\mu$.

**Lemma 1.13. (Theorem 3.2.[13]).** Let L be a complete lattice and $L(\mu, R)$ be an L-ring. Then intersection of arbitrary family of ideals of $\mu$ is an ideal of $\mu$.

## 2. PRIME RADICAL

The notion of prime ideal was introduced in our paper [15]. In the same paper we established several results pertaining to this notion. Here we formulate the concept of prime radical of an ideal in an L-ring.

***Definition 2.1 (Definition 2.1 [15]).*** *Let R be a commutative ring and Let* $L(\mu, R)$ *be an L-ring. An ideal* $\eta \neq \mu$ *of* $\mu$ *is said to be a prime ideal of* $\mu$ *if for all* $x, y \in R$ *, either*

$$\eta(xy) \wedge \mu(x) \wedge \mu(y) = \eta(x) \wedge \mu(y) \text{ or } \eta(xy) \wedge \mu(x) \wedge \mu(y) = \eta(y) \wedge \mu(x).$$

***Definition 2.2.*** *Let R be a commutative ring and Let L be a complete lattice. Let* $L(\mu, R)$ *be an L-ring and* $\eta$ *be an ideal of* $\mu$*. Write* $Q_\eta = \{\nu \mid \nu \text{ is a prime ideal of } \mu \text{ and } \eta \subseteq \nu\}$*. The prime radical of* $\eta$ *in* $\mu$ *(denoted by P($\eta$)) is defined by*

$$P(\eta) = \begin{cases} \bigcap \{\nu \mid \nu \in Q_\eta\} & \text{if } Q_\eta \neq \phi \\ \mu & \text{if } Q_\eta = \phi \end{cases}$$

By Lemma 1.13, $P(\eta)$ is an ideal of $\mu$. Clearly $\eta \subseteq P(\eta) \subseteq \mu$.

***Lemma 2.3 (Corollary 2.3 [18]).*** *Let Let R be a commutative ring and* $L(\mu, R)$ *be an L-ring. Let* $\eta$ *be an ideal of* $\mu$ *with* $\eta(0) < \mu(0)$*. Then the L-subset* $\xi: R \to L$ *defined by*

$$\xi(x) = \mu(x) \wedge \eta(0), \ \forall \ x \in R$$

*is a prime ideal of $\mu$ with $\eta \subseteq \xi$.*

**Theorem 2.4.** *Let R be a commutative ring and let L be a complete lattice. Let* $L(\mu, R)$ *be an L-ring and Let $\eta$ be an ideal of $\mu$. Then* $P(\eta)(0) = \eta(0)$.

**Proof.** Two cases arise :

Case (i) $\eta(0) = \mu(0)$. If $Q_\eta = \phi$, then $P(\eta) = \mu$. Hence $P(\eta)(0) = \mu(0) = \eta(0)$.

If $Q_\eta \neq \phi$, then $P(\eta)(0) = \left( \bigcap_{v \in Q_\eta} v \right)(0) = \bigwedge_{v \in Q_\eta} v(0) \geq \eta(0)$. Also

$$P(\eta)(0) = \left( \bigcap_{v \in Q_\eta} v \right)(0) \leq \mu(0) = \eta(0).$$

Thus $P(\eta)(0) = \eta(0)$.

Case (ii) $\eta(0) < \mu(0)$. By Lemma 2.3, the L-subset $\xi : R \to L$ defined by

$$\xi(x) = \mu(x) \wedge \eta(0) \ , \ \forall \ x \in R$$

is a prime ideal of $\mu$ and $\eta \subseteq \xi$. Thus $\xi \in Q_\eta$ . Now

$$P(\eta)(0) = \left( \bigcap_{v \in Q_\eta} v \right)(0) \leq \xi(0) = \mu(0) \wedge \eta(0) = \eta(0).$$

Also $P(\eta)(0) = \bigwedge_{v \in Q_\eta} v(0) \geq \eta(0)$. Thus $P(\eta)(0) = \eta(0)$. ∎

We also introduced the notion of semiprime ideal of an L-ring and semiprime radical of an ideal in [15]. We recall the following definitions and results.

**Definition 2.5 ([15]).** Let R be a commutative ring and let $L(\mu, R)$ be an L-ring. An ideal $\eta \neq \mu$ of $\mu$ is said to be a semiprime ideal of $\mu$ if

$$\eta(x^n) \wedge \mu(x) = \eta(x), \ \forall x \in R \ \& \ \forall n \in Z^+.$$

**Theorem 2.6 ([15]).** *Let R be a commutative ring and let $L(\mu, R)$ be an L-ring. Let $\eta$ be a prime ideal of $\mu$. Then $\eta$ is a semiprime ideal of $\mu$.*

**Definition 2.7 ([15]).** *Let R be a commutative ring and let L be a complete lattice. Let $L(\mu, R)$ be an L-ring and let $\eta$ be an ideal of $\mu$. The Radical of $\eta$, denoted by $\sqrt{\eta}$, is defined by*

$$\sqrt{\eta}(x) = \bigvee_{n \in Z^+} \left[ \eta(x^n) \wedge \mu(x) \right] \ , \ \forall x \in R.$$

Clearly $\eta \subseteq \sqrt{\eta} \subseteq \mu$.

**Definition 2.8 (Definition 2.1 [18]).** *Let R be a commutative ring and let L be a complete lattice. Let $L(\mu, R)$ be an L-ring and Let $\eta$ be an ideal of $\mu$. The semiprime radical of $\eta$ (denoted by $S(\eta)$) is defined by*

$$S(\eta) = \begin{cases} \bigcap \{ v \mid v \in T_\eta \} & \text{if } T_\eta \neq \phi \\ \mu & \text{if } T_\eta = \phi , \end{cases}$$

*where* $T_\eta = \{ v \mid v \text{ is a semiprime ideal of } \mu \text{ and } \eta \subseteq v \}$.

**Theorem 2.9 (Theorem 2.5 [18]).** *Let L be a complete lattice and $L(\mu, R)$ be an L-ring. Let $\eta$ be an ideals of $\mu$. Then $\sqrt{\eta} \subseteq S(\eta) \subseteq \mu$.*

**Theorem 2.10.** *Let R be a commutative ring and Let L be a complete lattice. Let $L(\mu, R)$ be an L-ring and Let $\eta$ and $\theta$ be ideals of $\mu$. Then the following assertions hold :*

*(i)* $\sqrt{\eta} \subseteq S(\eta) \subseteq P(\eta) \subseteq \mu$,

*(ii)* *if $\eta \subseteq \theta$, then $P(\eta) \subseteq P(\theta)$,*

*(iii)* $P(\eta \cap \theta) \subseteq P(\eta) \cap P(\theta)$.

**Proof.** We have $S(\eta) = \begin{cases} \cap\{\nu \mid \nu \in T_\eta\} & \text{if } T_\eta \neq \phi \\ \mu & \text{if } T_\eta = \phi \end{cases}$,

where $T_\eta = \{\nu \mid \nu \text{ is a semiprime ideal of } \mu \text{ and } \eta \subseteq \nu\}$.

(i) Since by Theorem 2.6 every prime ideal of $\mu$ is semiprime, we have $Q_\eta \subseteq T_\eta$. Thus $S(\eta) \subseteq P(\eta)$. By Theorem 2.9,n $\sqrt{\eta} \subseteq S(\eta)$. Hence

$$\sqrt{\eta} \subseteq S(\eta) \subseteq P(\eta) \subseteq \mu.$$

(ii) Let $\eta \subseteq \theta$. We show that $Q_\theta \subseteq Q_\eta$. Let $\nu \in Q_\theta$. Then $\nu$ is a prime ideal of $\mu$ and $\theta \subseteq \nu$. Now $\eta \subseteq \theta \subseteq \nu$. Thus $\nu \in Q_\eta$ and hence $Q_\theta \subseteq Q_\eta$. Therefore

$$P(\eta) = \cap\{\nu \mid \nu \in Q_\eta\} \subseteq \cap\{\nu \mid \nu \in Q_\theta\} = P(\theta).$$

(iii) Since $\eta \cap \theta \subseteq \eta$ and $\eta \cap \theta \subseteq \theta$, by (ii) we have $P(\eta \cap \theta) \subseteq P(\eta) \cap P(\theta)$. ∎

***Theorem 2.11 ([15]).*** *Let R be a commutative ring and Let L be a complete lattice. Let $L(\mu, R)$ be an L-ring. Then ideal $\eta$ of $\mu$ is a semiprime ideal of $\mu$ if and only if $\sqrt{\eta} = \eta$.*

***Theorem 2.12 (Theorem 2.6[17]).*** *Let R be a commutative ring and let L be a complete lattice. Let $L(\mu, R)$ be an L-ring. Then the intersection of an arbitrary family of semiprime ideals of $\mu$ is a semi prime ideal of $\mu$.*

***Theorem 2.13.*** *Let R be a commutative ring and let L be a complete lattice. Let $L(\mu, R)$ be an L-ring and let $\eta$ be an ideal of $\mu$. Then $P(P(\eta)) = P(\eta) = \sqrt{P(\eta)}$.*

**Proof.** In order to prove that $P(P(\eta)) = P(\eta)$, it is sufficient to show that $Q_{P(\eta)} = Q_\eta$. Since $\eta \subseteq P(\eta)$, we have $Q_{P(\eta)} \subseteq Q_\eta$. For the reverse inclusion, let $v \in Q_\eta$ Then $v$ is a prime ideal of $\mu$ such that $\eta \subseteq v$. Hence $P(\eta) = \cap \{\xi \mid \xi \in Q_\eta\} \subseteq v$. So that $v \in Q_{P(\eta)}$.

Next if $Q_\eta = \phi$, then $P(\eta) = \mu$ and hence $\sqrt{P(\eta)} = \sqrt{\mu} = \mu = P(\eta)$. And if $Q_\eta \neq \phi$, then $P(\eta) = \cap \{v \mid v \in Q_\eta\}$. Since every prime ideal of $\mu$ is a semiprime ideal, by Theorem 2.12, $P(\eta)$ is a semiprime ideal of $\mu$. Thus by Theorem 2.11, we have $\sqrt{P(\eta)} = P(\eta)$. ∎

**Theorem 2.14 (Theorem 2.12 [15]).** *Let R be a commutative ring and let L be a complete Heyting algebra. Let $L(\mu, R)$ be an L-ring and $\eta$ be an ideal of $\mu$. Then $\sqrt{\eta}$ is an ideal of $\mu$.*

**Theorem 2.15 (Theorem 2.13 [15]).** *Let R be a commutative ring and let L be a complete lattice. Let $L(\mu, R)$ be an L-ring and Let $\eta$ and $\theta$ be ideals of $\mu$. Then*

$$\eta \subseteq \theta \Rightarrow \sqrt{\eta} \subseteq \sqrt{\theta}.$$

**Theorem 2.16.** *Let R be a commutative ring. Let L be a complete Heyting algebra and $L(\mu, R)$ be an L-ring. Let $\eta$ be an ideal of $\mu$. Then $P(\sqrt{\eta}) = P(\eta) = \sqrt{P(\eta)}$.*

**Proof.** By Theorem 2.14, $\sqrt{\eta}$ is an ideal of $\mu$. Since $\eta \subseteq \sqrt{\eta}$, by Theorem 2.10 (ii), we have $P(\eta) \subseteq P(\sqrt{\eta})$. Also by Theorem 2.10 (i), $\sqrt{\eta} \subseteq P(\eta)$. Thus by Theorem 2.10 (ii) and Theorem 2.13, we have $P(\sqrt{\eta}) \subseteq P(P(\eta)) = P(\eta)$. Hence $P(\sqrt{\eta}) = P(\eta)$. ∎

**Theorem 2.17.** *Let R be a commutative ring and L be a complete Heyting algebra. Let $L(\mu, R)$ be an L-ring. Let $\eta$ and $\theta$ be ideals of $\mu$ with $\eta(0) = \theta(0)$. Then*

$$\sqrt{\eta} + \sqrt{\theta} \subseteq P(\sqrt{\eta} + \sqrt{\theta}) = P(\sqrt{\eta + \theta}) = P(\eta + \theta) \ .$$

**Proof.** By Lemma 1.11, $\eta + \theta$ is an ideal of $\mu$ and $\eta \subseteq \eta + \theta$, $\theta \subseteq \eta + \theta$. By Theorem 2.14, $\sqrt{\eta}$, $\sqrt{\theta}$ and $\sqrt{\eta + \theta}$ are ideals of $\mu$. Now $\sqrt{\eta}(0) = \eta(0) = \theta(0) = \sqrt{\theta}(0)$. Thus by Lemma 1.11, $\sqrt{\eta} + \sqrt{\theta}$ is an ideal of $\mu$. By Theorem 2.15, $\sqrt{\eta} \subseteq \sqrt{\eta + \theta}$ and $\sqrt{\theta} \subseteq \sqrt{\eta + \theta}$, therefore by Lemma 1.10 and Theorem 2.10 (i), we have

$$\sqrt{\eta} + \sqrt{\theta} \subseteq \sqrt{\eta + \theta} + \sqrt{\eta + \theta} = \sqrt{\eta + \theta} \subseteq P(\eta + \theta).$$

Thus by Theorem 2.10 (ii) and Theorem 2.13, we have $P(\sqrt{\eta} + \sqrt{\theta}) \subseteq P(P(\eta + \theta)) = P(\eta + \theta)$. Again since $\eta + \theta \subseteq \sqrt{\eta} + \sqrt{\theta}$, by Theorem 2.10, we have $P(\eta + \theta) \subseteq P(\sqrt{\eta} + \sqrt{\theta})$. Thus $P(\sqrt{\eta} + \sqrt{\theta}) = P(\eta + \theta)$. ■

The study of primary ideals of an L-ring was also initiated in [13]. We recall here some of the results which are required for the development of this work.

**Definition 2.18 (Definition 2.3 [16]).** *Let R be a commutative ring and let* $L(\mu, R)$ *be an L-ring. An ideal* $\eta \neq \mu$ *of* $\mu$ *is said to be primary ideal of* $\mu$ *if for all* $x, y \in R$*, we have either*

$$\eta(x) \wedge \mu(y) \geq \eta(xy) \wedge \mu(x) \wedge \mu(y) \qquad (1.1)$$

*or* $\qquad \eta(y) \wedge \mu(x) \geq \eta(xy) \wedge \mu(x) \wedge \mu(y) \qquad (1.2)$

*or* $\qquad \eta(x^n) \wedge \mu(x) \wedge \eta(y^m) \wedge \mu(y) \geq \eta(xy) \wedge \mu(x) \wedge \mu(y), \qquad (1.3)$

*for some integers m, n > 1.*

**Theorem 2.19 (Theorem 2.5 [16]).** *Let R be a commutative ring and let* $L(\mu, R)$ *be an L-ring Let* $\eta$ *be an ideal of* $\mu$ *with* $\eta \neq \mu$*. Then* $\eta$ *is a* primary ideal of $\mu$ if and only if for each non-empty level subset $\eta_t$, either $\eta_t = \mu_t$ or $\eta_t$ is primary ideal of $\mu_t$.

**Theorem 2.20 (Theorem 2.7 [16]).** Let R be a commutative ring and let L be a complete lattice. Let $L(\mu, R)$ be an L-ring and $\eta$ be a primary ideal of $\mu$ and has sup property. Then $\sqrt{\eta}$ is a prime ideal of $\mu$.

**Theorem 2.21.** Let R be a commutative ring. Let L be a complete lattice and $L(\mu, R)$ be an L-ring. Let $\eta$ be a primary ideal of $\mu$ and $\eta$ has sup property. Then $P(\eta) = \sqrt{\eta} = S(\eta)$.

**Proof.** By Theorem 2.10 (i), $\sqrt{\eta} \subseteq S(\eta) \subseteq P(\eta)$. Since $\eta$ is a primary ideal of $\mu$ and $\eta$ has sup property, by Theorem 2.20, $\sqrt{\eta}$ is a prime ideal of $\mu$. Thus $\sqrt{\eta} \in Q_\eta$. Hence $P(\eta) = \cap \{ \nu \mid \nu \in Q_\eta \} \subseteq \sqrt{\eta}$. Consequently $P(\eta) = \sqrt{\eta} = S(\eta)$. $\blacksquare$

**Corollary 2.22.** Let R be a commutative ring. Let L be a complete lattice and let $L(\mu, R)$ be an L-ring . Let $\eta$ be a primary ideal of $\mu$ and $\eta$ has sup property. Then the intersection of all prime ideals of $\mu$ containing the ideal $\eta$ is a prime ideal of $\mu$.

**Proof.** Since $P(\eta) = \sqrt{\eta}$ and $\sqrt{\eta}$ is a prime ideal of $\mu$, $P(\eta)$ is a prime ideal of $\mu$. $\blacksquare$

**Theorem 2.23.** Let R be a commutative ring. Let L be a complete lattice and let $L(\mu, R)$ be an L-ring. Let $\eta$ be a prime ideal of $\mu$. Then $P(\eta) = \eta = \sqrt{\eta} = S(\eta)$.

**Proof.** Obvious.

**Theorem 2.24.** *Let R be a commutative ring and let L be a complete lattice. Let L(μ,R) be an L-ring and let R be a commutative ring. et η be an ideal of μ. Then* $P(S(\eta)) = P(\eta) = S(P(\eta))$ .

**Proof.** By Theorem 2.10 (i), we have $\eta \subseteq S(\eta) \subseteq P(\eta)$. Thus by Theorem 2.10(ii) and Theorem 2.13, we have $P(\eta) \subseteq P(S(\eta)) \subseteq P(P(\eta)) = P(\eta)$. Therefore $P(S(\eta)) = P(\eta)$.

Again since $P(\eta)$ is an ideal, by Theorem 2.10(i) and Theorem 2.13, we have $P(\eta) \subseteq S(P(\eta)) \subseteq P(P(\eta)) = P(\eta)$. Thus $S(P(\eta)) = P(\eta)$. ∎

**Theorem 2.25.** *Let R be a commutative ring and let* L *be a complete lattice. Let* $L(\mu, R)$ *be an L-ring and let* η *and* θ *be ideals of* μ *with* $\eta(0) = \theta(0)$ . *Then*

$$P(\eta) + P(\theta) \subseteq P(P(\eta) + P(\theta)) = P(\eta + \theta).$$

**Proof.** By Lemma 1.11, $\eta + \theta$ and $P(\eta) + P(\theta)$ are ideals of μ and $\eta \subseteq \eta + \theta$, $\theta \subseteq \eta + \theta$. By Theorem 2.10 (ii) we have $P(\eta) \subseteq P(\eta + \theta)$ and $P(\theta) \subseteq P(\eta + \theta)$. Thus by Lemma 1.10, we have

$$P(\eta) + P(\theta) \subseteq P(\eta + \theta) + P(\eta + \theta) = P(\eta + \theta).$$

Hence by Theorem 2.10(ii) and Theorem 2.13, we have

$$P(P(\eta) + P(\theta)) \subseteq P(P(\eta + \theta)) = P(\eta + \theta).$$

Now, $\eta + \theta \subseteq P(\eta) + P(\theta)$, therefore by Theorem 2.10 (ii), we have

$$P(\eta + \theta) \subseteq P(P(\eta) + P(\theta)).$$

Hence $P(P(\eta) + P(\theta)) = P(\eta + \theta)$. ∎

**Corollary 2.26.** *Let* R *be a commutative ring and let L be a complete Heyting algebra. Let* $L(\mu, R)$ *be an L-ring. Let* η *and* θ *be ideals of* μ *with* $\eta(0) = \theta(0)$ . *Then* $\sqrt{\eta + \theta} \subseteq \sqrt{P(\eta) + P(\theta)} \subseteq P(\eta + \theta)$ .

**Proof**: By Theorem 2.14, $\sqrt{\eta + \theta}$ is an ideal of $\mu$. By the above theorem, $P(\eta) + P(\theta) \subseteq P(\eta + \theta)$, therefore by Theorem 2.15 and Theorem 2.13, we have

$$\sqrt{P(\eta) + P(\theta)} \subseteq \sqrt{P(\eta + \theta)} = P(\eta + \theta).$$

Since $\eta + \theta \subseteq P(\eta) + P(\theta)$, by Theorem 2.15, we have $\sqrt{\eta + \theta} \subseteq \sqrt{P(\eta) + P(\theta)}$.

## 3. PRIMARY DECOMPOSITION

**Definition 3.1.** *Let R be a commutative ring and let* $L(\mu, R)$ *be an L-ring. Let $\eta$ be an ideal of $\mu$. If there exist primary ideals $\eta_1, \eta_2, ..., \eta_m$ of $\mu$ such that $\eta = \bigcap\limits_{i=1}^{m} \eta_i$, then $\eta = \bigcap\limits_{i=1}^{m} \eta_i$ is called a primary decomposition (or primary representation) of $\eta$ in $\mu$. It is called irredundant or reduced if $\bigcap \{\eta_j \mid j = 1,2,...,m, \, j \neq i\} \nsubseteq \eta_i$ and $P(\eta_i)$ are all distinct, $i = 1, 2, ..., m$.*

**Lemma 3.2.** *Let J be a primary ideal of a commutative ring R and let I be a subring of R. Then either $J \cap I = I$ or $J \cap I$ is a primary ideal of I.*

**Lemma 3.3.** *Let R be a commutative ring. Let L be a chain and let $L(\mu, R)$ be an L-ring. Let $r, t_1 \in L$ with $r < t_1$. Let J be a primary ideal of the subring $\mu_r^{>}$ of R. Define an L-subset $\nu : R \to L$ by*

$$\nu(x) = \begin{cases} t_1 & \text{if } x \in J \\ r & \text{if } x \in R - J. \end{cases}$$

*Then $\nu \cap \mu$ is a primary ideal of $\mu$.*

**Proof.** Let $t \in L$. Then

$$(\nu \cap \mu)_t = \nu_t \cap \mu_t = \begin{cases} J \cap \mu_t & \text{if } r < t \leq t_1 \\ \mu_t & \text{if } t \leq r. \end{cases}$$

To show that $\nu \cap \mu$ is a primary ideal of $\mu$, by Theorem 2.19 it is sufficient to show that each non-empty level subset $(\nu \cap \mu)_t$ is either $\mu_t$ or is a primary ideal of $\mu_t$. Let $(\nu \cap \mu)_t$ be a non-empty level subset and $(\nu \cap \mu)_t \neq \mu_t$. Then $r < t \leq t_1$ and $(\nu \cap \mu)_t = J \cap \mu_t$. Since $r < t$, $\mu_t \subseteq \mu_r^>$. Also, since $L(\mu, R)$ is an L-ring, by Lemma 1.4 $\mu_t$ is a subring of $\mu_r^>$. Now by Lemma 3.2, $J \cap \mu_t$ is a primary ideal of $\mu_t$. Hence by Theorem 2.19, $\nu \cap \mu$ is a primary ideal of $\mu$.■

**Lemma 3.4.** *Let R be a commutative ring and let L be a chain. Let* $L(\mu, R)$ *be an L-ring and let* $\eta$ *be an ideal of* $\mu$. *Then each non-empty strong level subset* $\eta_t^>$ *of* $\eta$ *is an ideal of subring* $\mu_t^>$ *of R.*

**Theorem 3.5.** *Let R be a commutative ring and let L be a chain. Let* $L(\mu,R)$ *be an L-ring and let* $\eta$ *be an ideal of* $\mu$. *Let* $\text{Im}\,\eta$ *be finite and the ideal* $\eta_t^>$ *of subring* $\mu_t^>$ *has primary decompositions in* $\mu_t^>$ *for all* $t \in \text{Im}\,\eta - \sup_{x \in R}\{\eta(x)\}$. *Then* $\eta$ *has a primary decomposition in* $\mu$.

**Proof.** Let $\text{Im}\,\eta = \{t_1, t_2, \ldots, t_n\}$ with $t_1 > t_2 > \ldots > t_n$. For each i, $1 \leq i \leq n-1$, define L-subsets $\eta^{(i)} : R \to L$ as follows:

$$\eta^{(i)}(x) = \begin{cases} t_1 & \text{if } x \in \eta_{t_i} \\ t_{i+1} & \text{if } x \in R - \eta_{t_i} \end{cases}.$$

We show that $\bigcap_{i=1}^{n-1} \eta^{(i)} = \eta$. For $x \in \eta_{t_i}$, we have $t_i \leq \eta(x) \leq t_1 = \eta^{(i)}(x)$. For $x \in R - \eta_{t_i}$, we have $\eta(x) \leq t_{i+1} = \eta^{(i)}(x)$. Thus $\eta \subseteq \eta^{(i)}$, $\forall i = 1, 2, \ldots, n-1$. Hence $\eta \subseteq \bigcap_{i=1}^{n-1} \eta^{(i)}$. Clearly $\text{Im}\left(\bigcap_{i=1}^{n-1} \eta^{(i)}\right) = \{t_1, t_2, \ldots, t_n\} = \text{Im}\,\eta$. Now for each i, $1 \leq i \leq n-1$, we have

$$\eta_{t_k}^{(i)} = \begin{cases} \eta_{t_i} & \text{if } 1 \le k \le i \\ R & \text{if } i+1 \le k \le n \end{cases}$$

where $k = 1, 2, \ldots, n$. Thus $\eta_{t_k}^{(1)} = R$, $\eta_{t_k}^{(2)} = R$, ..., $\eta_{t_k}^{(k-1)} = R$, $\eta_{t_k}^{(k)} = \eta_{t_k}$, $\eta_{t_k}^{(k+1)} = \eta_{t_{k+1}}$, .... $\eta_{t_k}^{(n-1)} = \eta_{t_{n-1}}$. Also $\eta_{t_k} \subseteq \eta_{t_{k+1}} \subseteq \ldots \subseteq \eta_{t_{n-1}} \subseteq \eta_{t_n} = R$. Hence

$$\left( \bigcap_{i=1}^{n-1} \eta^{(i)} \right)_{t_k} = \bigcap_{i=1}^{n-1} \eta_{t_k}^{(i)} = \eta_{t_k}.$$

Thus by Lemma 1.8, we have $\bigcap_{i=1}^{n-1} \eta^{(i)} = \eta$. By Lemma 1.4, each non-empty strong level subset $\mu_t^>$ is a subring of R. By the hypothesis the ideals $\eta_{t_{i+1}}^>$ can be decomposed as an intersection of primary ideals of $\mu_{t_{i+1}}^>$ for all $i = 1, 2, \ldots n-1$. Now $\eta_{t_{i+1}}^> = \eta_{t_i}$. Thus each $\eta_{t_i}$ can be decomposed as an intersection of primary ideals of $\mu_{t_{i+1}}^>$.

Let $\eta_{t_i} = I_1^i \cap I_2^i \cap \ldots \cap I_{m_i}^i$, where $I_1^i, I_2^i, \ldots, I_{m_i}^i$ are primary ideals of $\mu_{t_{i+1}}^>$. For $j = 1, 2, \ldots, m_i$, consider L-subsets $\nu_j^i : R \to L$ defined by

$$\nu_j^i(x) = \begin{cases} t_1 & \text{if } x \in I_j^i \\ t_{i+1} & \text{if } x \in R - I_j^i \end{cases}.$$

Clearly $\bigcap_{j=1}^{m_i} \nu_j^i = \eta^{(i)}$. Therefore

$$\eta = \bigcap_{i=1}^{n-1} \eta^{(i)} = \bigcap_{i=1}^{n-1} \left( \bigcap_{j=1}^{m_i} \nu_j^i \right).$$

Hence $\eta = \eta \cap \mu = \bigcap_{i=1}^{n-1} \bigcap_{j=1}^{m_i} \left( \nu_j^i \cap \mu \right)$. By Lemma 3.3, $\nu_j^i \cap \mu$'s are primary ideals of $\mu$. Thus $\eta$ can be decomposed as an intersection of primary ideals of $\mu$. ∎

We recall here a lemma from [15] which is required to established the converse of Theorem 3.5.

**Lemma 3.6 (Lemma 2.4 [16]).** *Let R be a commutation ring. An ideal I of R is primary if and only if, whenever* $xy \in I$*, we have either* $x \in I$ *or* $y \in I$ *or* $(x^n \, \& \, y^m \in I)$*, for some integers m,n > 1.*

**Lemma 3.7.** *Let R be a commutative ring and let L be a chain. Let* $L(\mu,R)$ *be an L-ring and let* $\eta$ *be a primary ideal of* $\mu$*. Then for each non-empty strong level subset* $\eta_t^>$*, either* $\eta_t^> = \mu_t^>$ *or* $\eta_t^>$ *is a primary ideal of the subring* $\mu_t^>$*.*

**Proof.** Let $\eta_t^>$ be a non-empty strong level subset. By Lemma 3.4, $\eta_t^>$ is an ideal of subring $\mu_t^>$. Suppose $\eta_t^> \neq \mu_t^>$. To show that $\eta_t^>$ is a primary ideal of $\mu_t^>$, let $xy \in \eta_t^>$ , $x, y \in \mu_t^>$. Then $\eta(xy) > t$, $\mu(x) > t$, $\mu(y) > t$. Since L is a chain, $\eta(xy) \wedge \mu(x) \wedge \mu(y) > t$. Since $\eta$ is a primary ideal of $\mu$, by Definition 2.18, we have either

$$\eta(x) \wedge \mu(y) \geq \eta(xy) \wedge \mu(x) \wedge \mu(y) \tag{1.4}$$

or $\quad \eta(y) \wedge \mu(x) \geq \eta(xy) \wedge \mu(x) \wedge \mu(y) \tag{1.5}$

or $\quad \eta(x^n) \wedge \mu(x) \wedge \eta(y^m) \wedge \mu(y) \geq \eta(xy) \wedge \mu(x) \wedge \mu(y) \tag{1.6}$

for some integers $m, n > 1$.

If condition (1.4) holds then

$$\eta(x) \geq \eta(x) \wedge \mu(y) \geq \eta(xy) \wedge \mu(x) \wedge \mu(y) > t.$$

Thus $x \in \eta_t^>$. If (1.5) holds, then $y \in \eta_t^>$. In case condition (1.6) is valid, we have

$$\eta(x^n) \wedge \mu(x) \wedge \eta(y^m) \wedge \mu(y) \geq \eta(xy) \wedge \mu(x) \wedge \mu(y) \geq t$$

for some integers $m, n > 1$.

Thus $x^n, y^m \in \eta_t^>$. Therefore, by Lemma 1.25 $\eta_t^>$ is a primary ideal of $\mu_t^>$. $\blacksquare$

**Lemma 3.8**. *Let $L$ be a chain and $X$ be a non-empty set. Let $\eta_1, \eta_2, ..., \eta_m \in L^X$. Then $\left(\bigcap_{i=1}^{m} \eta_i\right)_t^> = \bigcap_{i=1}^{m} (\eta_i)_t^>, \ \forall \, t \in L$.*

**Theorem 3.9.** *Let $R$ be a commutative ring. Let $L$ be a chain and $L(\mu, R)$ be an $L$-ring. Let $\eta$ be an ideal of $\mu$ and has a primary decomposition in $\mu$. Then each non-empty strong level subset $\eta_t^> (\neq \mu_t^>)$ has a primary decomposition in the subring $\mu_t^>$ of R.*

**Proof.** Suppose $\eta = \eta_1 \cap \eta_2 \cap ... \cap \eta_k$ be a primary decomposition of $\eta$ in $\mu$. Then for each i, $\eta_i$ is a primary ideal of $\mu$. Let $\eta_t^>$ be a non-empty strong level subset. By Lemma 3.4, $\eta_t^>$ is an ideal of the subring $\mu_t^>$ of $R$. Now, by Lemma 3.8 $\eta_t^> = (\eta_1)_t^> \cap (\eta_2)_t^> \cap ... \cap (\eta_k)_t^>$. Since $\eta_i$ is a primary ideal of $\mu$ by the above lemma, either $(\eta_i)_t^> = \mu_t^>$ or $(\eta_i)_t^>$ is a primary ideal of $\mu_t^>$. If for all $i = 1, 2, ..., k$, $(\eta_i)_t^> = \mu_t^>$, then $\eta_t^> = \mu_t^>$, which is a contradiction. Thus there is at least one i such that $(\eta_i)_t^>$ is a primary ideal of $\mu_t^>$. Write $A = \left\{ i \mid (\eta_i)_t^> \neq \mu_t^> \right\}$. Then $\eta_t^> = \bigcap \left\{ (\eta_i)_t^> \mid i \in A \right\}$ is a primary decomposition of $\eta_t^>$ in the subring $\mu_t^>$ of R. $\blacksquare$

**Corollary 3.10.** *Let $R$ be a commutative ring and let $L$ be a chain. Let $L(\mu, R)$ be an $L$-ring and Let $\eta$ be an ideal of $\mu$ with $Im \eta$ finite. Then $\eta$ has primary decomposition in $\mu$ if the ideals $\eta_t^>$ of subring $\mu_t^>$ has primary* decompositions in $\mu_t^>$ for all $t \in \mathrm{Im}\,\eta - \sup\limits_{x \in R}\{\eta(x)\}$. Conversely if $\eta$ has a primary decomposition in $\mu$, then ideals $\eta_t^>(\neq \mu_t^>)$ has primary decomposition in the subring $\mu_t^>$ for all $t \in \mathrm{Im}\,\eta - \sup\limits_{x \in R}\{\eta(x)\}$.

**Lemma 3.11.** Let X be a non-empty set. Let $\eta_1, \eta_2, ..., \eta_m \in L^X$. Then

$$\left(\bigcap_{i=1}^m \eta_i\right)_t = \bigcap_{i=1}^m (\eta_i)_t, \forall\, t \in L$$

**Theorem 3.12.** Let R be a commutative ring and let L be a chain. Suppose for each L-ring $L(\mu, R)$, every ideal of $\mu$ has a primary decomposition in $\mu$. Then each ideal of any subring of R has a primary decomposition in the subring. In particular every ideal of R has a primary decomposition.

**Proof.** Let J be a subring of R and I be an ideal of J. We show that I has a primary decomposition in the subring J. Take $t_1, r \in L$ with $r < t_1$. Define L-subsets $\mu$ and $\eta$ from R to L by

$$\mu(x) = \begin{cases} t_1 & \text{if } x \in J \\ r & \text{if } x \in R - J \end{cases}$$

$$\eta(x) = \begin{cases} t_1 & \text{if } x \in I \\ r & \text{if } x \in R - I. \end{cases}$$

Then $L(\mu,R)$ is an L-ring and $\eta$ is an ideal of $\mu$. By the hypothesis $\eta$ has a primary decomposition in $\mu$. Let $\eta = \eta_1 \bigcap \eta_2 \bigcap ... \bigcap \eta_m$, where $\eta_1, \eta_2, ..., \eta_m$ are primary ideals of $\mu$ and $\eta_i \subsetneq \mu, \forall\, i = 1, 2, ..., m$. Then by Lemma 3.11

$$I = \eta_{t_1} = (\eta_1)_{t_1} \bigcap (\eta_2)_{t_1} \bigcap ... \bigcap (\eta_m)_{t_1}.$$

Now we show that $(\eta_1)_{t_1}, (\eta_2)_{t_2}, ...., (\eta_m)_{t_1}$ are primary ideals of $\mu_{t_1} = J$. Next we show that $(\eta_i)_{t_1} \neq \mu_{t_1}$, $\forall i = 1,2,..., m$. Suppose that $(\eta_i)_{t_1} = \mu_{t_1} = J$ for some i. Now $R = \eta_r \subseteq (\eta_i)_r \subseteq \mu_r = R$. Thus $(\eta_i)_r = R$. Hence $\mu \subseteq \eta_i \subseteq \mu$. Therefore $\eta_i = \mu$, which contradicts that $\eta_i \subsetneq \mu$. Thus $(\eta_i)_{t_1} \subsetneq \mu_{t_1}$, $\forall \ i = 1,2,..., m$. Since $\eta_i$ is a primary ideal of $\mu$, by Theorem 2.19 $(\eta_i)_{t_1}$ is a primary ideal of $\mu_{t_1}$ for all i=1,2,...,m. Hence I has a primary decomposition in J. ∎

**Lemma 3.13.** *Let R be a commutative ring .If an ideal J of R has a primary decomposition in R and I is a subring of R with* $J \subsetneq I$, *then J has a primary decomposition in I.*

**Proof.** Suppose $J = I_1 \cap I_2 \cap ... \cap I_m$ be a primary decomposition in R. Then $I_1, I_2, ..., I_m$ are primary ideals of R. Now

$$J = J \cap I = (I_1 \cap I) \cap (I_2 \cap I) \cap ... \cap (I_m \cap I).$$

By Lemma 3.2, either $I_i \cap I = I$ or $I_i \cap I$ is primary ideal of I. Now, if $I_i \cap I = I$ for all $i = 1,2,..., m$, then $J = I$, which is a contradiction. Thus, for at least one i, $I_i \cap I \neq I$. Write $B = \{i \mid I_i \cap I \neq I\}$. Then $J = \cap \{I_i \cap I \mid i \in B\}$ is a primary decomposition of J in I. ∎

**Corollary 3.14.** *Let R be a commutative ring and let L be a chain. Let* $L(\mu, R)$ *be an L-ring and let* $\eta$ *be an ideal of* $\mu$. *Let* $\mathrm{Im}\,\eta$ *be finite and ideals* $\eta_t \neq R$ *has primary decompositions in R for all* $t \in \mathrm{Im}\,\eta - \underset{x \in R}{\mathrm{Inf}} \{\eta(x)\}$ *and* $\eta_t^> \neq \mu_t^>$ *for all* $t \in \mathrm{Im}\,\eta - \underset{x \in R}{\sup}\{\eta(x)\}$. *Then* $\eta$ *has a primary decomposition in* $\mu$.

**Proof.** Suppose $\operatorname{Im}\eta = \{t_1, t_2, ..., t_n\}$ with $t_1 > t_2 > ... > t_n$. By the hypothesis $\eta_{t_i}$ has a primary decomposition in R for all $i = 1, 2, ..., n-1$. Now $\eta_{t_{i+1}}^{>} = \eta_{t_i}$, for all $i = 1, 2, ..., n-1$. By Lemma 1.4, $\mu_{t_{i+1}}^{>}$ is a subring of R for all $i = 1, 2, ..., n-1$. By Lemma 3.13, $\eta_{t_{i+1}}^{>}$ has a primary decomposition in $\mu_{t_{i+1}}^{>}$, for all $i = 1, 2, ..., n-1$. That is, $\eta_t^{>}$ has a primary decomposition in $\mu_t^{>}$ for all $t \in \operatorname{Im}\eta - \sup_{x \in R}\{\eta(x)\}$. By Theorem 3.5, $\eta$ has a primary decomposition in $\mu$. ∎

**Lemma 3.15 (Lemma 2.10 [15]).** *Let R be a commutative ring. Let L be a complete lattice and let* $L(\mu, R)$ *be an L-ring. Let* $\eta$ *be an ideal of* $\mu$ *and* $\eta$ *has sup property. Then* $\left(\sqrt{\eta}\right)_t = \sqrt{\eta_t} \cap \mu_t$, $\forall\ t \in L$.

**Theorem 3.16.** *Let R be a commutative ring and let L be a complete lattice. Let L($\mu$,R) be an L-ring. Let* $\eta$ *be an ideal of* $\mu$ *and has primary decomposition* $\eta = \eta_1 \cap \eta_2 \cap ... \cap \eta_k$ *in* $\mu$. *Then each non-empty level subset* $\eta_t (\neq \mu_t)$ *has a primary decomposition in the subring* $\mu_t$ *of R. Further if*

*(1)* *each* $\eta_i$ *has sup property,*

*(2)* *for some non-empty level subset* $\eta_t(\neq \mu_t)$, $\eta_t = (\eta_1)_t \cap (\eta_2)_t \cap ... \cap (\eta_k)_t$ *is a reduced primary decomposition of* $\eta_t$ *in* $\mu_t$,

*then* $\eta = \eta_1 \cap \eta_2 \cap ... \cap \eta_k$ *is a reduced primary decomposition of* $\eta$.

**Proof.** Since $\eta = \eta_1 \cap \eta_2 \cap ... \cap \eta_k$, by Lemma 3.11 we have $\eta_t = (\eta_1)_t \cap (\eta_2)_t \cap ... \cap (\eta_k)_t$. Since for each i, $\eta_i$ is a primary ideal of $\mu$, by Theorem 2.19, either $(\eta_i)_t = \mu_t$ or $(\eta_i)_t$ is a primary ideal of $\mu_t$. If for all i=1,2,...k, $(\eta_i)_t = \mu_t$, then $\eta_t = (\eta_1)_t \cap (\eta_2)_t \cap (\eta_k)_t = \mu_t$, which is a contradiction. Thus there exists at least one i such that $(\eta_i)_t$ is a primary ideal of $\mu_t$. Write $A = \{i | (\eta_i)_t \neq \mu_t\}$. Then $\eta_t = \bigcap \{(\eta_i)_t \mid i \in A\}$ is a primary representation of level subset $\eta_t$ in the subring $\mu_t$ of R.

In view of condition (3), suppose for some t, $\eta_t = (\eta_1)_t \bigcap (\eta_2)_t \bigcap ... \bigcap (\eta_k)_t$ is a reduced primary decomposition of the non-empty level subset $\eta_t$ in the subring $\mu_t$ of R. Then $(\eta_i)_t$ is a primary ideal of the subring $\mu_t$ for all $i = 1,2,...,k$. Now, we show that $\bigcap \{\eta_j \mid j = 1,2,..., k,\ j \neq i\} \nsubseteq \eta_i$ for all $i = 1,2,..., k$. If possible, $\bigcap \{\eta_j | j = 1,2,...., k,\ j \neq i\} \subseteq \eta_i$ for some i. Then

$$(\eta_i)_t \supseteq \bigcap \{(\eta_j)_t | j = 1,2,...,k,\ j \neq i\}.$$

This contradicts that $\eta_t = (\eta_1)_t \bigcap (\eta_2)_t \bigcap .... \bigcap (\eta_k)_t$ is a reduced primary decomposition of $\eta_t$ in $\mu_t$.

Next, we show that $P(\eta_i) \neq P(\eta_j)$ for $i \neq j$. Since $\eta_t = (\eta_1)_t \bigcap (\eta_2)_t \bigcap ... \bigcap (\eta_k)_t$ is a reduced primary decomposition of $\eta_t$ in the subring $\mu_t$ of R, we have

$\{x \in \mu_t \mid x^n \in (\eta_i)_t \text{ for some } n \in Z^+\} \neq \{x \in \mu_t | x^n \in (\eta_j)_t \text{ for some } n \in Z^+\}$ for $i \neq j$. Thus $\sqrt{(\eta_i)_t} \bigcap \mu_t \neq \sqrt{(\eta_j)_t} \bigcap \mu_t$ for $i \neq j$. Hence, by Lemma 3.15 $\left(\sqrt{\eta_i}\right)_t \neq \left(\sqrt{\eta_j}\right)_t$ for $i \neq j$. By Theorem 2.21, we have $P(\eta_i) = \sqrt{\eta_i}$ for all i = 1,2,...,k. Thus $P(\eta_i)_t = \left(\sqrt{\eta_i}\right)_t \neq \left(\sqrt{\eta_j}\right)_t = P(\eta_j)_t$ for $i \neq j$. Hence $P(\eta_i) \neq P(\eta_j)$ for $i \neq j$. Thus $\eta = \eta_1 \bigcap \eta_2 \bigcap ... \bigcap \eta_k$ is a reduced primary decomposition of $\eta$ in $\mu$.∎

***Corollary 3.17.*** *Let R be a commutative ring. Let L be a complete lattice and let L($\mu$,R) be an L-ring. Let* $t_0 = \underset{x \in R}{\wedge} \mu(x)$. *Let $\eta$ be an ideal of $\mu$ and has a primary decomposition*

$$\eta = \eta_1 \cap \eta_2 \cap \ldots \cap \eta_k.$$

*If $\eta_{t_0} \neq \mu_{t_0} = R$, then ideal $\eta_{t_0}$ of R has a primary decomposition in R. Further if*

*(1)    each $\eta_i$ has sup property,*

*(3)    $\eta_{t_0} = (\eta_1)_{t_0} \cap (\eta_2)_{t_0} \cap \ldots \cap (\eta_k)_{t_0}$ is a reduced primary decomposition of $\eta_{t_0}$, then $\eta = \eta_1 \cap \eta_2 \cap \ldots \cap \eta_k$ is a reduced primary decomposition of $\eta$ in $\mu$.*

**Corollary 3.18.** *Let R be a commutative ring. Let L be a complete lattice and let L($\mu$,R) be an L-ring. Let $\left\{ t \in L \mid t < \bigwedge_{x \in R} \mu(x) \right\}$ be non-empty. If every ideal of $\mu$ has a primary decomposition, then every ideal of R has a primary decomposition.*

**Proof.** Write $t_0 = \bigwedge_{x \in R} \mu(x)$. Then $\mu_{t_0} = R$. Choose $t_1 \in L$ such that $t_1 < t_0$. Let $I \neq R$ be an ideal of R. Define $\eta : R \to L$ by

$$\eta(x) = \begin{cases} t_0 & \text{if } x \in I \\ t_1 & \text{if } x \in R - I \end{cases}.$$

By Theorem 1.7, $\eta$ is an ideal of $\mu$ and $\eta_{t_0} = I \neq R$. By the hypothesis, $\eta$ has a primary decomposition in $\mu$. Thus by Theorem 3.16, $\eta_{t_0}$ has a primary decomposition in $\mu_{t_0} = R$. That is the ideal I of R has a primary decomposition.